\documentclass[a4paper,12pt]{article}
\usepackage{amsmath}
\usepackage{ulem}
\usepackage{calc}
\usepackage{vmargin}
\usepackage{amstext}
\usepackage{graphicx}
\usepackage{amsfonts}
\usepackage{amssymb}

\numberwithin{equation}{section}

\begin{document}

\begin{center}
{\Large {\bf Ideals of partial differential equations}}
\end{center}
\begin{center} \textbf{O.V. Kaptsov }\end{center}
\begin{center}
Institute of Computational Modelling, Siberian Branch of the Russian Academy of Sciences and Siberian Federal University
\end{center}

We propose a new algebraic approach to study compatibility of partial differential equations. The approach uses concepts from commutative algebra, algebraic geometry and  Gr\"{o}bner bases to clarify crucial notions concerning compatibility such as passivity (involution) and reducibility. One obtains sufficient conditions for a differential system to be passive and prove that such systems generate manifolds in the jet space. Some examples of constructions of passive systems associated with the $\sinh$-Gordon equation are given.

\section{Introduction}

 There are currently no methods to study the general systems of partial differential equations. Therefore it is necessary to investigate special classes of equations.  For example, the linear systems of homogeneous first order differential equations with one unknown function form one of  well-studied classes \cite{Goursat, Meleshko}. 

In the beginning of the twentieth century, French mathematicians Riquier, Janet, and Cartan made significant progress in studying a broad class of partial differential equations \cite{Riquier,Janet,Cartan}.  
Over the past several decades, new tools and terminology coming from differential geometry, differential and commutative algebra began to be applied in the formal theory of differential equations \cite{Kolchin,Pommaret,Vinograd}. 
It is now becoming increasingly important to consider algorithmic problems
of the theory of differential equations  \cite{Kondratieva,Seiler}. 
Some algorithms are implemented in computer algebra systems such as Maple, Reduce, Mathematica. 

In the papers \cite{Kaptsov1,Kaptsov2}, we used tools from the algebraic geometry and Gr\"{o}bner bases to study local properties of analytic partial differential equations. Here we consider smooth case. 
Some our notions can be explained by means of an example. 
Consider the $n+m$-dimensional real space $\mathbb{R}^{n+m}$ equipped with the natural coordinate functions $x_1,\dots,x_n,y_1,\dots,y_m$ and the standard topology.
Denote by $\mathfrak{F}(V)$ the algebra of smooth functions on an open set $V\subset\mathbb{R}^{n+m}$ and denote by $\mathfrak{F}_a$ the algebra of germs of smooth functions at a point $a\in\mathbb{R}^{n+m}$. 
A subset $B=\{f_1,\dots,f_m\}$ of $\mathfrak{F}(V)$ is called a normalized set, if each function $f_i\in B$ is of the form
$$f_i = y_i + g_i(x) ,$$
where the function $g_i$ can depends only on $x_1,\dots,x_n$. We say that an ideal of the algebra $\mathfrak{F}(V)$ is soft if it is generated by an normalized set. It is easy to give analogous definitions in the case of the algebra $\mathfrak{F}_a$.

The goal of this paper is to present an algebraic technique for studying compatibility of smooth partial differential equations.
Section 2 deals with  the infinite-dimensional space $\mathbb{R}^T$ of all the maps $T \longrightarrow \mathbb{R}$  equipped with the product topology (where $T$ is a countable set). 
To each open set $V$ of the space $\mathbb{R}^T$ one associates an algebra $\mathcal{F}(V)$ of smooth functions on $V$ depending only on finitely many variables. 
The set of all germs of these functions at a point $a\in\mathbb{R}^T$ forms a local algebra $\mathcal{F}_a$.
Next we define the appropriate normalized sets  and soft ideals in the algebras  $\mathcal{F}(V)$ and  $\mathcal{F}_a$. It turns out that every normalized set leads to a manifold in $a\in\mathbb{R}^T$.

We use the following notations: $\mathbb{N}_k =\{ 1,\dots,k \} $,
$\mathbb{N}$ is the set of all non-negative integers in Section 3 and consider
the infinite jet space $\mathbb{J}=\mathbb{R}^T$ with $T=\mathbb{N}_n \cup (\mathbb{N}_m\times \mathbb{N}^n) $. Then a system of partial differential equations is a subset of the algebra $\mathcal{F}(V)$.
We define passive systems of  partial differential equations at a point and on an open set in $\mathbb{J}$. These notations are analogous to  Gr\"{o}bner bases  \cite{Gathen}, but our definition does not apply any ranking.
We give simple properties of such systems.

In Section 4 we introduce the basic tool for study passive systems.
One of these is a stratified set which is given by a partition and a monoid acting on the set. Any stratified set must satisfy certain 
compatibility conditions. 
The monoid $(\mathbb{N}^n,+)$ acts on the algebras $\mathcal{F}(V)$ and  $\mathcal{F}_a$ by means of derivations. 
The stratification allows us to introduce  reductions of functions and germs modulo differential systems. We define reducibility conditions at
a point and on an open set in $\mathbb{J}$.

The crucial theorems are given in Section 5. We prove that if
a differential system $S$ is a normalize set and satisfies  
reducibility conditions at a point then it generates a soft ideal and it is passive. Furthermore, if the system  satisfies reducibility conditions on an open set then the orbit of $S$ leads to a manifold in the infinite jet space $\mathbb{J}$. 
At the end of our paper we give examples of  passive systems dealing with $\sinh$-Gordon equation. We briefly discuss how our construction can be applied to difference equations.

\section{Normalized sets in an algebra of smooth functions}

We shall use the following notations $\mathbb{R}$, for the set of all real numbers, $\mathbb{N}$, for the set of all non-negative integer,
$\mathbb{N}_k$, for the set $\{1,2,\dots,k\}$. 
Let $T$ be a denumerable set; the space of maps $z: T \longrightarrow
\mathbb{R}$ is denoted as $\mathbb{R}^T$ and equipped with the product topology. 
  This topology is obtained by taking as 
a neighborhood basis for a point $a\in \mathbb{R}^T$, sets of the form
\begin{equation} \label{U(a,rho)} 
U(a_{\tau},\rho) = \{z\in \mathbb{R}^{T}: |z_{t_i} - a_{t_i}|< \rho_i, \, i\in\mathbb{N}_k  \},
\end{equation}
where $t_i\in T$, $\rho_i > 0$, $\rho=(\rho_1,\dots,\rho_k)$,  $a_{\tau}=\{a_{t_1},\dots,a_{t_k} \}$ is a set of $k$ coordinates of the point $a$; $z_{t_1},\dots,z_{t_k}$ are $k$ coordinates of the point $z$.
 The functions $y_t : \mathbb{R}^T \longrightarrow \mathbb{R}$ defined by $y_t(z)= z(t)$, $t\in T$, are {\it the standard coordinate functions} (variables). 
 The set $Y =\{y_t \}_{t\in T}$ is {\it the standard coordinate system} on  $\mathbb{R}^T$. 

Let $V$ be an open set in $\mathbb{R}^T$ and let $\mathcal{F}(V)$ be the $\mathbb{R}$-algebra of real functions on $V$ that depend on finitely many variables and are smooth (i.e. they have derivatives of all orders) as functions of a finite number of variables. 
Suppose a function $f\in \mathcal{F}(V)$ depends on some set of 
variables, then this set denotes by $iv f$.
When $H$ is a subset of $\mathcal{F}(V)$, we shall use the  notation 
\begin{equation} \label{iv H}
 iv H =  \{iv f: f\in H \}.
\end{equation}

The family $\{\mathcal{F}(V)\}_{V\subset \mathbb{R}^T} $ gives rise to the sheaf  $\mathcal{F}$ of smooth functions on $\mathbb{R}^T$. 
For each point $a\in \mathbb{R}^T $ a stalk $\mathcal{F}_a$ of the sheaf is a $\mathbb{R}$-algebra of germs of smooth functions at $a$. 
Given a function $f\in \mathcal{F}(V)$, then its germ at $a$ is denoted 
$\tilde{f}_a$ or $\tilde{f}$ for simplicity. 

Each stalk $\mathcal{F}_a$ of the sheaf  $\mathcal{F}$ is a local algebra.
Indeed, if $\tilde{f}\in \mathcal{F}_a$ and $\tilde{f}(a) \neq 0$, then 
$1/\tilde{f}\in\mathcal{F}_a$ and $\tilde{f}$ does not belong to 
any proper ideal of the algebra. Hence the set  
\begin{equation} \label{M_a}
\mathfrak{M}_a = \{\tilde{f}\in \mathcal{F}_a: \tilde{f}(a)=0\} 
\end{equation}
is a unique maximal ideal of 
$\mathcal{F}_a$.

We shall say that the germ $\widetilde{f}\in \mathcal{F}_a$ depends on 
$\tilde{y}_t$ if there is a neighborhood $V$ of $a$ such that any representative $f$ of $\tilde{f}$ depends on $y_t$ in every  neighborhood $V^{\prime}\subset V$ of $a$. Assume a germ $\tilde{f}\in \mathcal{F}_a$ depends on a set of variables, then this set denotes by $iv \tilde{f}$.

\noindent
{\bf Definition 2.1.} A set $B\subset \mathcal{F}(V)$ is called {\it normalized} if the following conditions hold:

(i) any function $f\in B$ can be written $f=y_t+ g $, where the coordinate functions $y_t$ form a set $\mathcal{L}$ and the functions $g $ do not depend on elements of $\mathcal{L}$; 

(ii) if $f_1 = y_t + g_1, f_2 = y_t + g_2\in B$, then $f_1 =f_2$.
The elements of the sets $\mathcal{L}$ and $Y\setminus\mathcal{L}$ are called {\it principal} and {\it parametric} respectively.

We shall give a  similar definition for germs. Let $\tilde{Y_a}$ denote the set of germs of the coordinate functions at $a$.

\noindent
{\bf Definition 2.2.} A set $\tilde{B}\subset \mathcal{F}_a$ is called {\it normalized} if the following conditions hold:

(i) every germ $\tilde{f}\in \tilde{B}$ can be written $\tilde{f}=\tilde{y_t}+ \tilde{g} $, where the germs $\tilde{y_t}$ form a set $\tilde{\mathcal{L}}\subset \tilde{Y_a} $ and the germs $\tilde{g} $ do not depend on elements of $\tilde{\mathcal{L}}$;

(ii) if $\tilde{f_1} = \tilde{y_t} + \tilde{g_1}, \tilde{f_2} = \tilde{y_t} + \tilde{g_2}\in \tilde{B}$, then $\tilde{f_1} =\tilde{f_2}$.
The elements of the set $\tilde{\mathcal{L}}$ are called the {\it principal variables} and elements of the set $\tilde{Y}_a\setminus \tilde{\mathcal{L}}$ are {\it parametric variables}.

\noindent
{\bf Proposition 2.3.} Suppose $\tilde{f_i}= \tilde{y_{t_i}}+ \tilde{g_i}$, $i\in \mathbb{N}_k$, are some elements of a normalized set
 $\tilde{B}\subset \mathcal{F}_a$ and a germ $\tilde{F}\in \mathcal{F}_a$
 depends on $\tilde{y_{t_1}},\dots,\tilde{y_{t_k}}$. 
  Then there exist germs $\tilde{q_1},\dots,\tilde{q_k}\in \mathcal{F}_a$ and a unique germ  $\tilde{r}\in \mathcal{F}_a$ which does not depend on $\tilde{y_{t_1}},\dots,\tilde{y_{t_k}}$    such that
\begin{equation} \label{tF1}
\tilde{F} = \sum_{i=1}^{k} \tilde{q_i}\tilde{f_i} + \tilde{r} .
\end{equation}

{\it Proof.} Suppose the germs $\tilde{F},\tilde{f_1},\dots,\tilde{f_k}$ depend on $\tilde{y_{t_1}},\dots,\tilde{y_{t_n}}$. 
From the Mather division theorem \cite{BrLa}, we obtain  
$$ \tilde{F} = \tilde{q_1}\tilde{f_1} + \widetilde{r}_1 ,   $$
where $\widetilde{q}_1\in \mathcal{F}_a$; $\widetilde{r}_1\in \mathcal{F}_a$ does not depend on $\widetilde{y}_{t_1}$.
Using this theorem to the germ $\widetilde{r}_1$ yields 
$$ \tilde{F} = \tilde{q_1}\tilde{f_1} +\tilde{q_2}\tilde{f_2} + \widetilde{r}_2 ,   $$
where $\widetilde{r}_2$ does not depend on $\widetilde{y}_{t_1}, \widetilde{y}_{t_2}$. Continuing in the same way, we derive (\ref{tF1}).

One needs to verify uniqueness of $\tilde{r}$. Assume there exists another representation of $\tilde{F}$  
\begin{equation} \label{tF2}
\tilde{F} = \sum_{i=1}^{k} \tilde{q_i}^{\star}\tilde{f_i} + \tilde{r}^{\star}.
\end{equation}
It follows from (\ref{tF1}) and (\ref{tF2}) that
$$\tilde{f} = \sum_{i=1}^{k} \tilde{h_i}\tilde{f_i}            $$
with $\tilde{f} = \tilde{r}^{\star} - \tilde{r}$,
 $\tilde{h_i}=  \tilde{q_i} - \tilde{q_i}^{\star}$.
Let $f, h_i, f_i$ be representatives of the germs $\tilde{f}, \tilde{h_i},\tilde{f_i}$ . 
Then there is a neighborhood of $a$ such that 
\begin{equation} \label{f}
f = \sum_{i=1}^{k} h_if_i .
\end{equation}
Next we introduce new variables 
\begin{equation} \label{nv}
y^{\prime}_{t_1} = f_1, \dots,  y^{\prime}_{t_k} = f_k .   
\end{equation}
Since $f_i = y_{t_i} + g_i$ in some neighborhood of $a$, we can find $y_{t_i}$ from (\ref{nv}) and substitute  in  the expression (\ref{f}). 
Then we may write 
$$f = \sum_{i=1}^{k} \bar{h_i}y_{t_i}^{\prime} , $$
where $\bar{h_1},\dots,\bar{h_k}$ are some smooth functions while $f$ can only depend on $y_{t_{k+1}},\dots,y_{t_n}$.
 Assuming that 
$$y^{\prime}_{t_1} = 0, \dots, y^{\prime}_{t_k} = 0 ,$$
we have $f=0$ and therefore $\tilde{r}=\tilde{r}^{\star}$. 

\noindent
{\bf Proposition 2.4.} {\it Let $B\subset \mathcal{F}(V)$ be a normalized set. Assume that a function $F\in \mathcal{F}(V)$ is a polynomial in some principal variables $y_{t_1},\dots,y_{t_k}$ of $B$ with coefficients depending only on parametric variables.   
 Then there is  a unique function $r\in \mathcal{F}(V)$ not depending on the principal variables and some functions $q_1,\dots,q_k\in\mathcal{F}(V)$
such that 
\begin{equation} \label{F=}
F= \sum q_i f_i +r ,
\end{equation}
where $f_i=y_{t_i}+g_i\in B$. }

\noindent
{\it Proof.} The function $F$ is a polynomial in the principal variables $y_{t_1},\dots,y_{t_k}$ and the functions $f_1,\dots,f_k$ are polynomials of the first degree with coefficients $1$. 
Then we can obtain (\ref{F=}) using the multivariate division
with remainder \cite{Gathen}, although $\mathcal{F}(V)$ is not a field.
Moreover, the function $r$ does not depend on the principal variables and lies in $\mathcal{F}(V)$.

The uniqueness of $r$ can be proved as in  Proposition 2.3.  Suppose that the function $F$ is written in the other form
\begin{equation} \label{F=2}
F= \sum q^{\prime}_i f_i +r^{\prime} ,
\end{equation}  
where the function $r^{\prime}$ does not depend on the principal variables.
Then from (\ref{F=}) and (\ref{F=2}) we have
\begin{equation} \label{r=2}
r^{\prime\prime}= \sum q^{\prime\prime}_i f_i    
\end{equation}
with $r^{\prime\prime}=r^{\prime}-r$ and $q^{\prime\prime}_i=q^{\prime}_i - q_i$. 
 Under the transformation
 \begin{equation} \label{cv}
 y^{\prime}_{t_1} = f_1, \dots,  y^{\prime}_{t_k} = f_k .   
 \end{equation}
 the relation (\ref{r=2}) becomes
$$ r^{\prime\prime}= \sum q^{\star}_i y^{\prime}_{t_i}  ,      $$
where $ q^{\star}_1,\dots,q^{\star}_k\in\mathcal{F}(V)$, while the function $r^{\prime\prime}$ does not depend on $y_{t_1},\dots,y_{t_k}$.
Setting 
$$y^{\prime}_{t_1} = 0, \dots, y^{\prime}_{t_k} = 0 ,$$
we obtain $r^{\prime\prime}=r^{\prime}-r=0$.  

\noindent
{\it Remark.} Inserting the values $y_{t_1}=-g_1,\dots,y_{t_k}=-g_k$
 in the function $F$, we obtain the function $r$. We say that $r$ is a {\it normal form} of $F$ modulo $B$. 

A general definition of a smooth manifold is in \cite{Zharin}, but we shall only consider embedded submanifolds of $\mathbb{R}^T$.

\noindent
{\bf Definition 2.5} Let $V$ be an open set in $\mathbb{R}^T$.
A map $\phi: V \rightarrow \mathbb{R}^T$ is called smooth on $V$ if 
for all $t\in T$ the functions  $\phi_t = y_t\circ\phi$ are smooth on $V$. 
 
Let $V$, $V^{\prime}$ be open sets in $\mathbb{R}^T$. We say that a map $\psi: V \rightarrow V^\prime$ is a diffeomorphism 
if $\psi$ carries $V$ homeomorphically onto $V^{\prime}$ and if $\psi$ and
$\psi^{-1}$ are smooth. 
 If $T^{\star}\subset T$, then a set
$$ C_{T^{\star}} = \{z\in \mathbb{R}^T: z(t)=0, \forall t\in T \}  $$
is called a coordinate subspace of $\mathbb{R}^T$.
We shall assume that any subset $Q\subset \mathbb{R}^T$ 
is equipped with a topology induced from that of $\mathbb{R}^T$.


\noindent
{\bf Definition 2.6.} A subset $M\subset \mathbb{R}^T$ is called a {\it smooth manifold} if for any $a\in M$ there are a neighborhood $V\subset \mathbb{R}^T$, an open subset $V^\prime \subset \mathbb{R}^T$,
 and a diffeomorphism $\phi: V \rightarrow V^\prime$ such that 
$$ \phi(V\cap M) = V^\prime \cap C_{T^{\star}} , $$
where  $C_{T^{\star}}$ is a coordinate subspace of $\mathbb{R}^T$.

\noindent
{\bf Proposition 2.7.}  Assume that $\{g_t\}_{t\in T^{\prime}}$ is a family of smooth function on an open subset $W\subset \mathbb{R}^{T^{\prime\prime}}$ 
with  $T^{\prime\prime} = T \setminus T^{\prime}$  and denote by $V$ the open set $W\times\mathbb{R}^{T^\prime}$ in $\mathbb{R}^T$. Then a set
$B = \{y_t + g_t\}_{t\in T^{\prime}}\subset \mathcal{F}(V) $ is normalized and the set 
$$Z(B) = \{z\in V: f(z)=0, f\in B \}     $$
is a manifold in $\mathbb{R}^T$.

{Proof.}  Let $ \phi: V\rightarrow \mathbb{R}^T$ be a map given by
$$ y^\prime_t = y_t + g_t , \ y^\prime_s = y_s  $$
with $t\in  T^{\prime}, s\in T^{\prime\prime}$. Then the inverse map is of the form
$$  y_t = y^\prime_t - g_t , \ y_s = y^\prime_s   . $$
It is easy to see that  
$$ \phi(V\cap Z(B)) = V \cap \mathbb{R}^{T^{\prime\prime}} , $$
and hence $Z(B)$ is a manifold.

\section{Passive differential systems}

We now introduce the notion of a soft ideal which plays an important role in the study of compatibility of equations.

\noindent
{\bf Definition 3.1.}

(i) We say that a proper ideal $I$ of an algebra $\mathcal{F}(V)$ is {\it soft},
 if there is a normalized set $B\subset \mathcal{F}(V)$ to generate  the ideal.  The set $B$ is called  a {\it normalized system of generators} of $I$.
 
(ii) Let $J$ be a proper ideal of an algebra $\mathcal{F}_a$.
A normalized subset $\tilde{B}\subset \mathcal{F}_a$  generating the ideal $J$ is called a {\it normalized system of generators}
  of $J$ and we say that the ideal is {\it soft}.
  
We recall that a derivation in an algebra $A$ over $\mathbb{R}$ is a map 
$\mathcal{D}: A \longrightarrow A   $  such that 
$$ \mathcal{D}(ab) = a\mathcal{D}(b) + \mathcal{D}(a)b, \quad
\mathcal{D}(k_1a+k_2b)=k_1\mathcal{D}(a)+ k_2\mathcal{D}(b)  $$
for all $a,b\in A$ and for all $k_1,k_2\in\mathbb{R}$.

\noindent
The next proposition describes an arbitrary derivation of the algebra of germs $\mathcal{F}_a$.

\noindent
{\bf Proposition 3.2.} 
Let $ \mathcal{D},  \bar{\mathcal{D}}$ be derivations of the algebra 
$\mathcal{F}_a$ such that $\mathcal{D}(y_t)=\bar{\mathcal{D}}(y_t)$ for all
$y_t\in Y$. Then $\mathcal{D}=\bar{\mathcal{D}} $ and
\begin{equation} \label{Df}
\mathcal{D}(\tilde{f}) = \sum_{t\in T} \frac{\partial\tilde{f}}{\partial\tilde{y_t}}\mathcal{D}(\tilde{y_t}) , \quad \forall \tilde{f}\in \mathcal{F}_a .
\end{equation}

{\it Proof}. Repeating the proof Theorem 4.2 (a variant of Hadamard's lemma) in \cite{BrLa}, we see that  the set $\tilde{Y}$ of germs $\{\tilde{y_t}\}_{t\in T}$ at $a\in \mathbb{R}^T$ generates the maximal ideal (\ref{M_a}). It follows from the Proposition 8.16
\cite{Jac} that $\mathcal{D}=\bar{\mathcal{D}}$.  It is easy to see that the expression (\ref{Df}) gives the derivations of $\mathcal{F}_a$.
 Even though the formula (\ref{Df}) involves an infinity summation, when applying $\mathcal{D}$ to any germ $\tilde{f}$, only finitely many terms are need.

Now we proceed to consider differential equations. Further, assume that
$$ T = \mathbb{N}_n \cup (\mathbb{M}\times \mathbb{N}^n) ,    $$
where $\mathbb{M}=\mathbb{N}_m$ or $\mathbb{M}=\mathbb{N}$. By $\mathbb{J}$ denote the space $\mathbb{R}^T$,  and call it the {\it jet space}.
The canonical coordinate functions on $\mathbb{J}$ are denoted by $x_1,\dots,x_n$, $u^i_\alpha$, where $i\in \mathbb{M}, \alpha\in\mathbb{N}^n$. The Cartesian coordinate system $Y$ on $\mathbb{J}$
is decomposed into two sets    
\begin{equation} \label{XU}
X=\{x_1,\dots,x_n\} , \qquad U=\{u^i_\alpha\}^{i\in\mathbb{M}}_{\alpha\in\mathbb{N}^n} . 
\end{equation}

The elements $e_1 =(1,0,\dots,0), \dots, e_n =(0,\dots,1)$ are generators of the monoid $\mathbb{N}^n$. Introduce derivations $\mathcal{D}_1,\dots,\mathcal{D}_n$ on the algebras $\mathcal{F}(V)$, $\mathcal{F}_a$ so that
\begin{equation} \label{D_i}
\mathcal{D}_i f = \frac{\partial f}{\partial x_i}  + 
 \sum_{j\in\mathbb{M}, \alpha\in\mathbb{N}^n } \frac{\partial f}{\partial u^i_\alpha} u^i_{\alpha + e_j} , \quad
\mathcal{D}_i \tilde{f} = \frac{\partial\tilde{f}}{\partial \tilde{x_i}}  + 
 \sum_{j\in\mathbb{M}, \alpha\in\mathbb{N}^n } \frac{\partial \tilde{ f}}{\partial \tilde{u^i_\alpha}} \tilde{u^i}_{\alpha + e_j} .
\end{equation}
Thus $\mathcal{F}(V)$ and $\mathcal{F}_a$ became differential algebras. 

We, following Ritt's lead \cite{Ritt}, call the coordinate functions $u^i_0$ the {\it indeterminates}, and $u^i_\alpha$, where $\alpha\neq 0$, the partial derivatives of $u^i_0$.

\noindent
{\bf Definition 3.3.}
We shall say that a subset $S\subset\mathcal{F}(V)$ is a {\it differential system} on an open set $V\subset \mathbb{J}$, if any function $f\in S$ depends on at least one of the  partial derivatives.
If $\mathbb{M}=\mathbb{N}_m$ then we say that $S$ is a {\it system with finite number of indeterminates}, but if  $\mathbb{M}=\mathbb{N}$ then we get a {\it system in infinitely many indeterminates}.

Let $W$ be an open set in $\mathbb{R}^n$, 
$h : W \rightarrow \mathbb{R}^{\mathbb{M}}$ a smooth map  with components $h_m$ for $m\in \mathbb{M}$. 
Then a map $h^{\infty} :  W \rightarrow \mathbb{J}$
 whose components are $x_i,h^\alpha_m=\mathcal{D}^\alpha (h_m)$ for $i\in\mathbb{N}_n$, $m\in \mathbb{M}$, $\alpha \in \mathbb{N}^n$
 is called the {\it infinite prolongation graph} of  $h$.

\noindent
{\bf Definition 3.4.} Let $S$ be a differential system on an open subset
 $V\subset \mathbb{J}$. A smooth map  $h: W \rightarrow \mathbb{R}^{\mathbb{M}}$ is called a solution  of a differential system $S$ if the following conditions hold:
$$ (1)\quad h^{\infty}(W)\subset V ,  \qquad  (2)\quad f\circ h^{\infty} = 0 , \quad \forall f \in S.$$

\noindent
{\it Remarks.} In other words, the map $h$ is a solution of the system $S$ if under substitution of $\mathcal{D}^\alpha (h_i)$ for $u^i_\alpha$ every function $f\in S$ vanishes. A germ of a solution is defined in the obvious way.

 An ideal of the algebra $\mathcal{F}(V)$ generated by a set 
$\{\mathcal{D}^\alpha(f): f\in S ,  \alpha\in\mathbb{N}^n \} $
we shall denote by $\langle\langle S\rangle\rangle$. 
Similarly, let $\tilde{S}_a$ be a set of germs of functions in $S\subset\mathcal{F}(V)$ at $a$.  An ideal of the algebra $\mathcal{F}_a$ generated by the set
$\{\mathcal{D}^\alpha(\tilde{f}): \tilde{f}\in \tilde{S} ,  \alpha\in\mathbb{N}^n \} , $
denoted by $\langle\langle \tilde{S} \rangle\rangle_a$. 

It is obvious that a map $h$ is a solution of a differential system $S$ if and only if $f\circ h^{\infty} = 0$ for all $f\in \langle\langle S\rangle\rangle$. There are some cases in which it is convenient to deal  with other differential system $S^\prime$ such that 
$\langle\langle S^\prime \rangle\rangle = \langle\langle S \rangle\rangle$.
In particular, such examples arise when we consider compatible systems of differential equations of the first order for a single unknown function \cite{Meleshko}.

Recall that if $G$ and $H$ are sets, then $G$ acts on $H$ in case there is a mapping $\psi: G\times H \longrightarrow H$. The mapping $\psi$ is called a action. When $\psi$ is fixed, then $gh$ denotes $\psi(g,h)$. 
The monoid $(\mathbb{N}^n,+,0)$ acts on the algebras $\mathcal{F}(V)$, 
$\mathcal{F}_a$ by 
$$ \alpha f = \mathcal{D}^\alpha f , \qquad  \alpha \tilde{f} = \mathcal{D}^\alpha \tilde{f} , \qquad \forall\alpha\in\mathbb{N}^n
\forall f \in\mathcal{F}(V) \forall \tilde{f} \in\mathcal{F}_a . $$
The sets 
$$ O(f)= \{\mathcal{D}^\alpha f: \alpha \in \mathbb{N}^n \}  , \quad O(\tilde{f})=  \{\mathcal{D}^\alpha \tilde{f}: \alpha \in \mathbb{N}^n\}                         $$
are orbits of a function $f$ and a germ $\tilde{f}$ under $(\mathbb{N}^n,+,0)$.

\noindent
{\bf Definition 3.5.}

(i) A germ $\tilde{f}\in \mathcal{F}_a$ of the form $\tilde{f} =\tilde{u}^i_\alpha + \tilde{g} $ is called {\it solvable} with respect to
$\tilde{u}^i_\alpha$, if the germ $\tilde{g}$ does not depend on elements of the orbit $O(\tilde{u}^i_\alpha)$. 

(ii) A function $f=u^i_\alpha\in\mathcal{F}(V)$ is {\it solvable} with respect to $u^i_\alpha$, if the function $g$ does not depend on elements of the orbit $O(u^i_\alpha)$.

Suppose a germ $\tilde{f}\in \mathcal{F}_a$ is solvable with respect to
$\tilde{u}^i_\alpha$. Then the germ $\tilde{u}^i_\alpha$ is denoted by 
$st\tilde{f}$. Let $\tilde{S}_a$ be a set of  solvable germs at a point $a$, then we shall use the notation $st\tilde{S} = \{st\tilde{f}: \tilde{f}\in \tilde{S}_a\} $. The same notation is used for functions.

 
  
\noindent
{\bf Definition 3.6.} A differential system $S\subset\mathcal{F}(V)$ is called {\it passive at $a\in V$}, if the ideal $\langle\langle \tilde{S} \rangle\rangle_a$ is smooth, the set $\tilde{S}_a$ consists of  solvable germs, and a set of principal variables of a normalized system of the ideal $\langle\langle \tilde{S} \rangle\rangle_a$ coincides with the orbit $O(st \tilde{S}_a) $.
The system $S$ is {\it passive on $V$} if every function in $S$ is solvable,
the ideal  $I=\langle\langle S \rangle\rangle$ is smooth, and a set of principal variables of a normalized system of the ideal $I$ coincides with the orbit $O(st S)$.





\noindent
{\bf Proposition 3.7.} Let $S\subset\mathcal{F}(V)$ be a passive system at $a\in V$. Then every non-zero germ $\tilde{f}\in \langle\langle \tilde{S} \rangle\rangle_a$ depends on at least one element of $O(st \tilde{S}) $.

{\it Proof.} 
Suppose there exists a non-zero germ $\tilde{f}\in \langle\langle \tilde{S} \rangle\rangle_a$ which does not depend on elements of $O(st \tilde{S}) $. 
If $\tilde{B}$ is the normalized system of generators of the ideal 
$\langle\langle \tilde{S} \rangle\rangle_a$, then we have
\begin{equation} \label{1f=}
\tilde{f} = \sum_{i=1}^{p} \tilde{q}_i \tilde{b}_i ,
\end{equation}
where $\tilde{q}_i\in \mathcal{F}_a$, $\tilde{b}_i\in \tilde{B}$.

It follows from (\ref{1f=}) and the definition 3.5  that there are a neighborhood $W$ of the point $a$ and  representatives of the germs
$\tilde{f}, \tilde{q}_i,  \tilde{b}_i$ such that 
\begin{equation} \label{2f=}
f = \sum_{i=1}^{p} q_i b_i ,
\end{equation}
in any neighborhood $V^\prime\subset W$ of $a$. It is cleat that the functions $b_1,\dots,b_p$ form a normalized set. 
 Suppose the functions $q_1,\dots,q_p,b_1,\dots,b_p$ depend on some variables  $y_1,\dots,y_r\in Y$.
  and assume $y_1=st\ b_1,\dots,y_p=st\ b_p$ for simplicity.
We make a change of variables
\begin{equation} \label{change}
y_1^\prime = b_1 , \dots, y_p^\prime = b_p.
\end{equation}
Since $b_i = y_i + g_i$ in some neighborhood of $a$, we can find $y_i$ from (\ref{change}) and substitute  in  the expression (\ref{2f=}). As a result we obtain
\begin{equation} \label{3f=}
f = \sum_{i=1}^{p} \bar{q}_i y^\prime_i ,
\end{equation}
where $f$ can only depend on $y_{p+1}, \dots, y_r$, and $\bar{q}_i$, $i=1,\dots,p$, are  functions in $\mathcal{F}_a$.
Substituting $y^\prime_i=0$ into (\ref{3f=}) yields  $f=0$, and this contradicts the assumption that the germ $\tilde{f}$ is non-zero.

\section{Stratified sets and reductions}

We need a convenient criterion for recognizing passive systems.
For this purpose, we shall introduce additional tools.
Recall that a preorder $\preceq$ is a binary relation that is reflexive and transitive. A strict partial order $\prec$ is a binary relation that is irreflexive, transitive.

In what follows, we shall deal with a well-ordered set $\Gamma$.
Every partition $\{H_\gamma\}_{\gamma\in \Gamma}$  of a set $H$
gives rise to a preorder and a strict partial order on $H$ as follows:
\begin{equation} \label{prec1}
h_1\preceq h_2 \quad \Longleftrightarrow \quad \exists\gamma_1,\gamma_2\in\Gamma
( \gamma_1\leq \gamma_2 \wedge h_1\in H_1 \wedge h_2\in H_2  ) ,  
\end{equation}
\begin{equation} \label{prec2}
h_1\prec h_2 \quad \Longleftrightarrow \quad \exists\gamma_1,\gamma_2\in\Gamma
( \gamma_1 < \gamma_2 \wedge h_1\in H_1 \wedge h_2\in H_2  ) . 
\end{equation}
In this case we say that the set $H$ is equipped with a induced strict partial order. We also say that a monoid $G$ acts on the set $H$ if there
exists a map $(g, x) \longrightarrow gx$ of $G\times H$ into $H$ satisfying 
$$eh = h, \qquad (g_1g_2)h= g_1(g_2h) \quad \forall h\in H\forall g_1,g_2\in G,$$
where $e$ is the identity of $G$.

\noindent
{\bf Definition 4.1.} Suppose $\{H_\gamma\}_{\gamma\in \Gamma}$ is a partition of a set $H$ equipped with a induced strict partial order, $G$ is a monoid acting on $H$. We shall say that $H$
is a {\it stratified $G$-set} if for all $g\in G$ it satisfies the following conditions : 

$$ 1) \quad \forall\gamma \forall h_1 \forall h_2 \exists\gamma^\prime (h_1,h_2\in H_\gamma \Longrightarrow gh_1,gh_2\in H_{\gamma^\prime})  ;   $$
$$ 2) \quad h_1 \prec h_2 \Longrightarrow gh_1\prec gh_2 ;    $$
$$ 3) \quad h\prec gh \qquad \forall h\in H \forall g\in G \quad(g\neq e),  $$
where $e$ is the identity of $G$. 

\noindent {\it Remark.} The above definition is a generalization of ranking \cite{Kolchin}.

Define an action of the monoid $(\mathbb{N}^n,+)$ on the set of coordinate function $U$ by the rule
$$  \beta u^i_\alpha = u^i_{\alpha+\beta}   \qquad 
\forall \alpha,\beta\in\mathbb{N}^n \forall i\in\mathbb{M}     $$
with $\mathbb{M}=\mathbb{N}_m$ or $\mathbb{M}=\mathbb{N}$.
 It is easy to see that $U = \bigcup_{n\in\mathbb{N}} U_n $ with
$U_n = \{u^i_\alpha\in U: |\alpha|=n \} $
 gives an example of stratified $\mathbb{N}^n$-set.  

Let $V$ be an open set in $\mathbb{J}$ and $X$ is given by (\ref{XU}).
We consider two sets 
\begin{equation} \label{F_X}
\mathcal{F}(V)_X=  \{f\in\mathcal{F}(V): iv f \subset X\}, \qquad  
\hat{\mathcal{F}}(V) = \mathcal{F}(V)\setminus \mathcal{F}(V)_X .
\end{equation}
We shall indicate how a partition $\{U_\gamma \}_{\gamma\in \Gamma}$ of the set $U$ leads to a partition of $\hat{\mathcal{F}}(V)$.

Consider sets
\begin{equation} \label{Ygamma}
Y_\gamma = X\cup (\bigcup_{\gamma_0\leq\gamma^\prime\leq\gamma } U_{\gamma^\prime} ) , \qquad \gamma_0=min\{\gamma\in\Gamma \}
\end{equation}
which form an ascending chain of subsets of $Y$. 
The sets 
\begin{equation}\label{J_gamma}
J^{\gamma}= \{z\in \mathbb{J}: y(z)=0, \forall y\in (Y\setminus Y_\gamma) \}  ,  
\end{equation}
\begin{equation}\label{Fgamma}
\mathcal{F}^\gamma(V) = \{f\in \mathcal{F}(V) : iv(f)\subset Y_{\gamma}  \}     
\end{equation}
also form ascending chains of subspaces and subalgebras respectively.
This chain of subalgebras generates a partition $\{\Phi^\gamma(V) \}_{\gamma\in\Gamma}$ of the set $\hat{\mathcal{F}}(V)$ , where  
\begin{equation}\label{Phigamma}
\Phi^\gamma(V) = \mathcal{F}^\gamma(V) \setminus (\bigcup_{\gamma_0<\gamma^\prime <  \gamma} \mathcal{F}^{\gamma^\prime}(V) ) , \qquad \Phi^{\gamma_0}(V) = \mathcal{F}^{\gamma_0}(V) \setminus  
\mathcal{F}(V)_X .
\end{equation}
Let us consider three set of germs
\begin{equation}\label{Fgamma_a}
\mathcal{F}^\gamma_a = \{\tilde{f}\in \mathcal{F}_a : iv(\tilde{f})\subset \tilde{Y}_{\gamma} \}  ,    
\end{equation}
\begin{equation} \label{F_Xa}
\mathcal{F}_{aX} =  \{\tilde {f}\in\mathcal{F}_a: iv \tilde{f} \subset \tilde{X}\}, \qquad  
\hat{\mathcal{F}}_a = \mathcal{F}_a\setminus \mathcal{F}_{aX} .
\end{equation}
A partition $\{\Phi^\gamma_a\}_{\gamma\in \Gamma}$ of the set $\hat{\mathcal{F}}_a$ is given by
\begin{equation}\label{Phigamma_a}
\Phi^\gamma_a = \mathcal{F}^\gamma_a \setminus (\bigcup_{\gamma_0<\gamma^\prime < \gamma} \mathcal{F}^{\gamma^\prime}_a ) , \qquad \Phi^{\gamma_0}_a = \mathcal{F}^{\gamma_0}_a \setminus \mathcal{F}_{aX} .
\end{equation}

\noindent
{\bf Lemma 4.2.} {\it Suppose that the set $U$ (\ref{XU}) is a stratified $\mathbb{N}^n$-set. Then the sets $\hat{\mathcal{F}}(V)$ and $\hat{\mathcal{F}}_a$ are also stratified $\mathbb{N}^n$-sets.}

\noindent
{\it Proof.}
It suffices to check three requirements of a stratified set for generators of the monoid $\mathbb{N}^n$.
At first, we consider the set $\hat{\mathcal{F}}(V)$.
 To prove first property of a stratified set it will suffice to show 
the following statement. 
If $f_1,f_2\in \Phi^{\gamma}(V)$, then there exists an element $\gamma^\prime\in\Gamma$ such that functions $\mathcal{D}_k(f_1),\mathcal{D}_k(f_2)$ given by (\ref{D_i}) lie in $\Phi^{\gamma^\prime}(V)$. 

We remark that if $\frac{\partial f}{\partial u^i_\alpha}$ vanishes  
on some open set $\Omega$ in $\mathbb{J}$ then the function $f$ does not depend on $ u^i_\alpha$ in  $\Omega$.

Since $f_1,f_2\in \Phi^{\gamma}(V)$, then there are variables $u^i_\alpha, u^j_\beta\in U^\gamma$, and points $a_1, a_2\in V$ such that 
$$\frac{\partial f_1}{\partial u^i_\alpha}(a_1)\neq 0, \qquad 
 \frac{\partial f_2}{\partial u^j_\beta}(a_2)\neq 0  .     $$
It follows from assumption of our Lemma that for all $u^i_\alpha, u^j_\beta \in U^\gamma$ there exists $\gamma^\prime\in\Gamma$ such that 
$\mathcal{D}_k u^i_\alpha, \mathcal{D}_k u^j_\beta$ lie in
 $U^{\gamma^\prime}$. Thus,  we clearly obtain
 $ \frac{\partial f_1}{\partial u^i_\alpha}u^i_{\alpha+e_k} ,  
 \frac{\partial f_2}{\partial u^j_\beta} u^j_{\beta+e_k}\in \Phi^{\gamma^\prime}(V) $ and furthermore, 
  $\mathcal{D}_k f_1, \mathcal{D}_k f_2\in\Phi^{\gamma^\prime}(V) $.
In a similar manner, one can prove two other properties. 

We shall now prove that $\hat{\mathcal{F}}_a$ is also a stratified $\mathbb{N}^n$-set. At first,  we show that if $\tilde{f}\in \Phi^\gamma_a$
then for any representative $f$ of the germ $\tilde{f}$ there exists a 
an neighborhood $V^\star$ of $a$ such that for every  
neighborhood $V^\prime\subset V^\star$  of $a$ there are a variable $u^i_\alpha\in U^\gamma$ and a point $b\in V^\prime$ with 
$\frac{\partial f}{\partial u^i_\alpha }(b)\neq 0$. Suppose this is not the case. Then there exists a representative $f$ of the germ $\tilde{f}$
such that for every neighborhood $V^\star$ of $a$ there is a
neighborhood $V^\prime\subset V^\star$  of $a$ in which $\frac{\partial f}{\partial u^i_\alpha }(b) = 0$, for any variable $u^i_\alpha\in U^\gamma$ and every point $b\in V^\prime$. Therefore, the function $f$ does not depend
on variables $u^i\alpha\in U^\gamma$ in neighborhood $V^\prime$.
We thus get a contradiction to $\tilde{f}\in \Phi^\gamma_a$.
This implies that  $\frac{\partial\tilde{f}}{\partial \tilde{u}^i_\alpha } \neq 0$.

Let us prove the first property of a stratified set for $\hat{\mathcal{F}}_a$. Suppose that $\tilde{f}_1$ and $\tilde{f}_2$ lie in $\Phi^\gamma_a$. 
It suffices to show that $\mathcal{D}_k\tilde{f}_1$ and $\mathcal{D}_k\tilde{f}_2$ lie in $\Phi^{\gamma^\prime}_a$ for some $\gamma^\prime\in\Gamma$.
From assumption of this lemma, there is an element $\gamma^\prime\in\Gamma$ such that  
$$\mathcal{D}_k\tilde{u}^i_\alpha =\tilde{u}^i_{\alpha+e_k}, \quad
\mathcal{D}_k\tilde{u}^j_\beta =\tilde{u}^j_{\beta+e_k} \quad 
\forall \tilde{u}^i_\alpha, \tilde{u}^j_\beta \in \tilde{U}^\gamma .$$
It follows as above that there are variables $u^i_\alpha, u^j_\beta\in U^\gamma$, an element $\gamma^\prime\in\Gamma$ and a number $k\in\mathbb{N}$
such that germs
$\frac{\partial\tilde{f}_1}{\partial \tilde{u}^i_\alpha }\tilde{u}^i_{\alpha+e_k},
\frac{\partial\tilde{f_2}}{\partial \tilde{u}^j_\beta}\tilde{u}^j_{\beta+e_k}$ lie in $\Phi^{\gamma^\prime}_a$.
Hence $\mathcal{D}_k\tilde{f}_1,\mathcal{D}_k\tilde{f}_2\in \Phi^{\gamma^\prime}_a$.
The other properties are proved in the same vein. 

In what follows we shall suppose that $\hat{\mathcal{F}}(V)$ and  $\hat{\mathcal{F}}_a$ are  stratified $\mathbb{N}^n$-sets equipped with a induced strict partial order.  

\noindent
{\bf Definition 4.3.}

(i) A function $f= u^i_\alpha + g\in\mathcal{F}(V)$ is called {\it orderly solvable} (with respect to $u^i_\alpha$), if $g\prec u^i_\alpha$. The variable $u^i_\alpha$ is denoted by $lt f$ and is called {\it leading term of $f$}. 

(ii) A germ $\tilde{f}= \tilde{u}^i_\alpha + \tilde{g}\in\mathcal{F}_a$ is called {\it orderly solvable} (with respect to $\tilde{u}^i_\alpha$), if $\tilde{g}\prec \tilde{u}^i_\alpha$.  The germ $\tilde{u}^i_\alpha$ is  denoted by $lt \tilde{f}$ and is called {\it leading term of $\tilde{f}$}. 


It is clear that if a germ $\tilde{f}\in\mathcal{F}_a$ is  orderly solvable with respect to $\tilde{u}^i_\alpha$ then it is solvable with respect to $\tilde{u}^i_\alpha$ in terms of the definition 3.4. In the future, we suppose that $st \tilde{f}= lt \tilde{f}$ in this case. Furthermore, we assume that $st\tilde{S}=lt\tilde{S}$ for every orderly solvable local system
with  $lt \tilde{S}=\{lt f: f\in S\}$.

\noindent
{\bf Proposition 4.4.} {\it Let $\tilde{F}\in\mathcal{F}_a$ be a germ depending on 
$\tilde{u}^i_\beta$. Suppose that $\tilde{f}=\tilde{u}^i_\alpha+\tilde{g}$ is 
a orderly solvable germ with respect to $\tilde{u}^i_\alpha$
 and there exists $\delta\in\mathbb{N}^n$ satisfying $\beta=\alpha+\delta$. Then
there exists a unique germ $\tilde{r}\in\mathcal{F}_a$ and a germ $q\in\mathcal{F}_a$
such that 
\begin{equation}\label{F==} 
 \tilde{F}= \tilde{q}D^\delta\tilde{f}+ \tilde{r} , \quad 
 \tilde{u}^i_\alpha\notin iv\tilde{r} 
\end{equation} 
\begin{equation}\label{r=} 
\tilde{q}\preceq \tilde{F} , \quad \tilde{r} \preceq \tilde{F} .     
\end{equation} }
{\it Proof.} The germ $D^\delta \tilde{f}$ is equal to 
$\tilde{u}^i_\beta+D^\delta\tilde{g}$, where $D^\delta\tilde{g}\prec \tilde{u}^i_\beta$. Then from the Mather theorem \cite{BrLa}, we obtain 
(\ref{F==}). The uniqueness $\tilde{r}$ is proved  just as in the second part of Proposition 2.3. 
It is clear that 
$$  iv(\tilde{q})\subseteq (iv(\tilde{F})\cup iv(D^\delta\tilde {g} )  ) ,
\quad   iv(\tilde{r})\subset  (iv(\tilde{F})\cup iv(D^\delta\tilde {g})), \quad \tilde{u^i_\alpha} \notin iv\tilde{r} .
$$
Since $lt(D^\delta\tilde {g})= \tilde{u}^i_\alpha $, it follows that 
$D^\delta\tilde {g}\preceq \tilde{F} .$  The last relations lead to (\ref{r=}).

If the assumptions of Proposition 4.4 are satisfied, then we say that the {\it the germ $\tilde{F}$ reduces to $\tilde{r}$ modulo $\tilde{f}$} at $a$, 
denoted by $\tilde{F} \xrightarrow[\tilde{f}]{} \tilde{r}$. 

\noindent
{\bf Proposition 4.5.} {\it Let $F$ be a polynomial in $u^i_\beta$
with coefficients that do not depend on $u^i_\beta$ and lie in $\mathcal{F}(V)$. 
Assume that  $f=u^i_\alpha+g$ is a  orderly solvable function with respect to $u^i_\alpha$ and $\delta$ is a element in  $\mathbb{N}^n$ satisfying $\beta=\alpha+\delta$. Then    
there exists a unique function $\tilde{r}\in\mathcal{F}(V)$ and a function $q\in\mathcal{F}(V)$ such that 
\begin{equation}\label{2F==} 
 F= qD^\delta f + r , \quad 
 u^i_\alpha\notin iv r 
\end{equation} 
\begin{equation}\label{2r=} 
q\preceq F , \quad r \preceq F .     
\end{equation} }

\noindent
{\it Proof.} The relation (\ref{2F==}) follows from Proposition 2.4.
The inequalities (\ref{2r=}) are proved 
just as in the second part of Proposition 4.4.

If the assumptions of Proposition 4.5 are satisfied, then we say that the {\it the function $F$ reduces to the function $r$ modulo $f$} on $V$, 
denoted by $F \xrightarrow[f]{} r$. 

\noindent
{\bf Definition 4.6.} A differential system $S\subset\mathcal{F}(V)$ is called {\it weakly solvable}, if every function $f\in S$ is orderly solvable. We write
$lt S = \{lt f: f\in S\}$.

\noindent
{\bf Definition 4.7.}
Let $S\subset\mathcal{F}(V)$ be a weakly solvable differential system.

(i) We shall say that {\it a germ $\tilde{F}\in\mathcal{F}_a$  reduces to a germ $\tilde{r}\in\mathcal{F}_a$ modulo $\tilde{S}_a$}, written
   $\tilde{F} \xrightarrow[\tilde{S}]{} \left. \tilde{r}\right|_a$, if 
 there exists a consequence of germs $\tilde{r}_1,\dots,\tilde{r}_{k-1}\in\mathcal{\tilde{F}}_a$ such that
$$\tilde{F}\xrightarrow[\tilde{f}_1]{} \tilde{r}_1 \xrightarrow[\tilde{f}_2]{}\cdots \xrightarrow[\tilde{f}_{k-1}]{} \tilde{r}_{k-1}\xrightarrow[\tilde{f}_k]{} \tilde{r}  $$
with $\tilde{f}_1,\dots,\tilde{f}_k\in \tilde{S}_a$.

(ii) Let $S$ be a normalized set in $\mathcal{F}(V)$.
Suppose that $F\in\mathcal{F}(V)$ is a polynomial in $O(lt S)$ with coefficients depending only on variables in $O(Y\setminus ltS)$. We say that 
{\it $F$ reduces to a function $r\in\mathcal{F}(V)$ modulo $S$}, written
   $F \xrightarrow[S]{} r$, 
if there exists a consequence of functions $r_1,\dots,\tilde{r}_{k-1}\in\mathcal{F}(V)$ such that
$$F \xrightarrow[f_1]{} r_1 \xrightarrow[f_2]{}\cdots \xrightarrow[f_{k-1}]{} r_{k-1}\xrightarrow[f_k]{} r  $$
with $f_1,\dots,f_k\in S$.

Let us define a binary operation $\diamond$ on $\mathbb{N}^n$ by
$$\alpha\diamond\beta = (\mu_1,\dots,\mu_n) ,$$
where $\alpha=(\alpha_1,\dots,\alpha_n )$, $\beta =(\beta_1,\dots,\beta_n )$,
$\mu_i = max(\alpha_i, \beta_i) - \alpha_i$.
Suppose that functions $f_1,f_2\in\mathcal{F}(V)$ are orderly solvable with respect to $u^i_\alpha, u^i_\beta$  respectively and $\tilde{f_1}, \tilde{f_2}$
are their germs at $a\in V$. Then we define two differences
\begin{equation}\label{tau}
\tau(f_1,f_2)=D^{\alpha \diamond \beta} f_1 - D^{\beta\diamond \alpha }f_2,
\qquad \tau(\tilde{f_1},\tilde{f_2})=D^{\alpha \diamond \beta}\tilde{f_1} - D^{\beta\diamond \alpha }\tilde{f_2} .
\end{equation}

\noindent
{\bf Definition 4.8.} Let $S\subset\mathcal{F}(V)$ be a weakly solvable differential system.  

(i)  The system $S$ satisfies {\it  reducibility conditions at $a\in V$}, if 
\begin{equation}\label{tau1}
\tau(\tilde{f_1},\tilde{f_2})  \xrightarrow[\tilde{S}]{} \left. \tilde{0} \right|_a
\end{equation}
for each pair of functions $f_1,f_2\in S$ such that $lt f_1=u^i_\alpha$, 
$lt f_2 = u^i_\beta$.  

(ii) Let $S$ be a normalized set in $\mathcal{F}(V)$.
 We say that $S$ satisfies {\it reducibility conditions on $V$}, if 
\begin{equation}\label{tau2}
\tau(f_1,f_2)  \xrightarrow[S]{} 0
\end{equation}
for each pair of functions $f_1,f_2\in S$ such that $lt f_1=u^i_\alpha$, 
$lt f_2 = u^i_\beta$.

 Denote by $\mathbb{D}$ an algebra of operators such that 
every element of $\mathbb{D}$ can be written as a finite sum 
\begin{equation}\label{mathbbD}
P = \sum a_\alpha\mathcal{D}^\alpha
\end{equation}
with  $a_\alpha\in\mathbb{R}$. Let $RU$ be a vector space over $\mathbb{R}$
consisting of finite sums 
\begin{equation}\label{RU}
s= \sum b_i^\beta u^i_\beta , \qquad b_i^\beta\in\mathbb{R}.
\end{equation}
Define an action of $\mathbb{D}$ on $RU$ by letting  
$$Pu^i_\beta = \sum a_\alpha u^i_{\alpha+\beta} ,$$
and extending $P$ to $RU$ by linearity. 

\noindent
{\bf Definition 4.9.} Let $\bf{y}$ be an $k$-tuple $(y_{t_1},\dots,y_{t_k})$ of variables $y_{t_i}\in U$. An $k$-tuple
 ${\bf d}=(d_1,\dots,d_k)$ of operators in $\mathbb{D}$ is called {\it syzygy} of $\bf{y}$,  if
 $$  d_1 y_{t_1}+\cdots+  d_k y_{t_k} = 0 .$$ 
\noindent
It is clear that the syzygies of the $k$-tuple $\bf{y}$ constitute  a $\mathbb{D}$-module denoted by $Syz\ \bf{y}$. 

Suppose ${\bf y}=(y_{t_1},\dots,y_{t_k})\in U^k$ with $y_{t_i}=u^l_\alpha$ and  
$y_{t_j}=u^l_\beta$, then 
\begin{equation}\label{sigma}
\sigma_{ij} = \mathcal{D}^{\alpha\diamond\beta}e_i - 
\mathcal{D}^{\beta\diamond\alpha}e_j
\end{equation}
is a syzygy of $\bf{y}$. It is easy to show (see \cite{Kaptsov1}) that 
the syzygies (\ref{sigma}) generate  the $\mathbb{D}$-module $Syz\ \bf{y}$,
if a number of the indeterminates $u^i_0\in U$ is finite.

\noindent
{\bf Example.} Assume $m=1$ and $n=2$, so that $U=\{u_{(i,j)}: i,j\in\mathbb{N}\}$; take ${\bf y}=(u_{(0,1)},u_{(0,2)},u_{(1,1)} )$.
It is obvious that $(\mathcal{D}_2,-1,0)$, $(\mathcal{D}_1,0,-1)$ and 
$(0,\mathcal{D}_1,-\mathcal{D}_2)$ are syzygies of the $3$-tuple $\bf{y}$.

\section{Passivity criterion of differential systems}

In this section we give a sufficient condition for a differential system to be passive. Furthermore, we prove that a passive system generates a manifold in the jet space.

Let $S\subset\mathcal{F}(V)$ be weakly solvable differential system. We call a point $a\in\mathbb{J}$ equivalent to a point $b\in\mathbb{J}$, written $a\sim b$, if  $y(a)=y(b)$  for all coordinate functions $y\in Y\setminus O(lt S)$. 

\noindent
{\bf Theorem 5.1.} {\it Let $S=\{f_1,\dots,f_k\}\subset\mathcal{F}(V)$ be a differential system  with finite number of indeterminates. Suppose that $S$ is a normalized set and satisfies reducibility conditions (\ref{tau1}) at $a\in V$.
Then  the following properties hold:

(1) there is a unique point $b\sim a$ such that 
\begin{equation}\label{Db}  
D^\alpha f(b) = 0 , \quad \forall f\in S \ \forall \alpha\in\mathbb{N}^n ;
\end{equation}     


(2) the system $S$ is passive at any point $c\sim a$.} 

\noindent
{\it Proof.} Since $S$ is a normalized set, we conclude that the orbit $O(S)$ is a weakly solvable differential system.  
This gives rise to the uniqueness of the point $b$ satisfying  the  condition (\ref{Db}). 

We have shown above that a partition $\{U_\gamma \}_{\gamma\in \Gamma}$ of the set $U$ provides the ascending chain of subspaces $J^\gamma$ (\ref{J_gamma}), the chains of subalgebras $\mathcal{F}^\gamma(V)$  (\ref{Fgamma}),   $\mathcal{F}^\gamma_z$ (\ref{Fgamma_a}) and leads to the partitions $\{\Phi^\gamma(V)\}_{\gamma\in\Gamma}$ (\ref{Phigamma}),  $\{\Phi^\gamma_z\}_{\gamma\in\Gamma}$ (\ref{Phigamma_a}) with $z\in V$. We also recall that $Y_\gamma$ is defined by (\ref{Ygamma}).
Consider  linear mappings $\pi_\gamma : \mathbb{J}\longrightarrow J^\gamma$, where the coordinates of $\pi_\gamma(z)$ are given by
$$ y(\pi_\gamma(z))= y(z)  \quad \forall y\in Y_\gamma ; \qquad
 y(\pi_\gamma(z))=0 \quad \forall y\in Y\setminus Y_\gamma \ .$$
Recall that  $\tilde{S}_z$ is a set of germs of functions in $S$ at the point $z$.
We shall use the following notion:
$$\gamma_0 = min\{\gamma\in\Gamma: O(\tilde{S}_a)\cap\Phi^\gamma_a\neq 0\}, $$
$$O^\gamma_z=O(\tilde{S}_z)\cap\mathcal{F}^\gamma_z , \quad 
C^\gamma_z = O(\tilde{S}_z)\cap\Phi^\gamma_z .$$
It is obvious that
\begin{equation}\label{Ostar} 
O^{\gamma_\star}_z = C^{\gamma_\star}_z \cup 
(\bigcup_{\gamma_0\leq\gamma <\gamma_\star} C^\gamma_z )
\end{equation} 
for any $\gamma_\star>\gamma_0$.  
Let $\langle O^\gamma_z \rangle$ be an ideal of the algebra  $\mathcal{F}^\gamma_z$ generated by $O^\gamma_z$. 
  
We shall use transfinite induction to prove that for all $\gamma\geq\gamma_0$
the following properties hold:

(i) there exists a point $b_\gamma\sim \pi_\gamma(a)$ such that $\tilde{f}(b_\gamma)=0$
for all $\tilde{f}\in O^\gamma_{b_\gamma} $;

(ii) there exists a normalized system $\tilde{B}^\gamma_c$ of generators of the ideal 
 $\langle O^\gamma_c \rangle$ for any point $c\sim \pi_\gamma(a)$.

Assume that $\gamma=\gamma_0$ then two cases arise:
 
 1. All leading terms of germs in $C^{\gamma_0}_a$ are distinct.
 
 2. There exist at least two germs $\tilde{f}_i, \tilde{f}_j\in C^{\gamma_0}_a$ such that $lt \tilde{f}_i = lt \tilde{f}_j$.
 
It is clear that in the first case there is a point $b_{\gamma_0}\sim \pi_{\gamma_0}(a)$ such that  $\tilde{f}(b_{\gamma_0})=0$
for all $\tilde{f}\in C^{\gamma_0}_{b_{\gamma_0}}$. 
Furthermore,  the properties (ii) and (iii) are satisfied because $\tilde{B}^{\gamma_0}_c=\tilde{C}^{\gamma_0}_c$ is a normalized system of generators of the ideal  $\langle O^{\gamma_0}_c \rangle$ and $B^{\gamma_0}=S^{\gamma_0}$ is a normalized system of generators of the ideal $\langle O^{\gamma_0}(S) \rangle$. 
In the second case, there must be germs $\tilde{f}_i, \tilde{f}_j\in \tilde{C}^{\gamma_0}_a$ such that $lt \tilde{f}_i = lt \tilde{f}_j$. Then 
$\tilde{f}_i - \tilde{f}_j\in \mathcal{F}^{\gamma^\prime}_a$, where $\gamma^\prime< \gamma_0$, and
 $\tilde{f}_i - \tilde{f}_j  \xrightarrow[\tilde{S_a}]{} \tilde{0} $ according to the conditions of our theorem. Since $O^{\gamma^\prime}_a$ is the empty set then we have $\tilde{f}_i = \tilde{f}_j$.

 Assume that our statement is true for all $\gamma$ with $\gamma_0\leq\gamma < \gamma_\star$ and prove its for $\gamma = \gamma_\star$.
As above, we need to distinguish two cases:

1. All leading terms of germs in $C^{\gamma_\star}_a$ are distinct.

2. There exist two germs $\tilde{f}, \tilde{g}\in C^{\gamma_\star}_a$ such that $lt \tilde{f} = lt \tilde{g}$.

\noindent 
In the first case, the property (i) is trivially satisfied. According to the assumption of induction and the formula (\ref{Ostar}), the set 
$$G^{\gamma_\star}_c = C^{\gamma_\star}_c \cup 
(\bigcup_{\gamma_0\leq\gamma <\gamma_\star} \tilde{B}^\gamma_c )     $$
is a system of generators (not necessarily normalized) of the ideal $\langle O^{\gamma_\star}_c \rangle$ for any point $c\sim \pi_{\gamma_\star}(a)$.  
 
In the second case, there are two germs $\tilde{f}, \tilde{g}\in C^{\gamma_\star}_a$ with  $lt \tilde{f} = lt \tilde{g}$.  Then there exist  two germs $\tilde{f_p}, \tilde{f_q}\in \tilde{S}_a$ such that
$$ lt \tilde{f}= lt D^\mu \tilde{f_p}= lt \tilde{g} =lt D^\eta  \tilde{f_q} ,$$ 
where $D^\mu=D_1^{\mu_1}\cdots D_n^{\mu_n}$ and $D^\eta = D_1^{\eta_1}\cdots D_n^{\eta_n}$ are some differential monomials. Therefore, we have
\begin{equation}\label{ltfp} 
D^\mu (lt \tilde{f_p}) =  D^\eta(lt \tilde{f_q}) .
\end{equation} 
Denote by $\bf{y}$ an n-tuple constructed from all elements of the set
$lt \tilde{S}_a$. Assume that the elements $lt \tilde{f_p}$ and $lt \tilde{f_q}$ are the $i$-th and $j$-th  items in $\bf{y}$. It follows from 
(\ref{ltfp}) that $d = D^\mu e_i - D^\eta e_j$ is a syzygy of $\bf{y}$.
It is easy to see that there is differential monomial $D^\nu$ such that 
$d=D^\nu\sigma_{ij} $, where $\sigma_{ij}$ is one of the sygyzies (\ref{sigma}) generating $\mathbb{D}$-module $Syz\ \bf{y}$.

The difference $\tilde{f}-\tilde{g}$  reduces to the zero germ  modulo $\tilde{S}_a$. Indeed, the system $S$ satisfies reducibility conditions at $a$ by assumption, then we have
$$ D^\nu \sigma_{ij}(\tilde{f_p},\tilde{f_q})=  
\tilde{f}- \tilde{g}  \xrightarrow[\tilde{S_a}]{} \tilde{0} . $$
Next, any one of the germs $\tilde{f}, \tilde{g}$ is included in a set $gen^{\gamma_\star}_a$ while the other is not. 
In the same way we search for all pairs of germs in $C^{\gamma_\star}_a$ with equal leading terms, form the set $gen^{\gamma_\star}_a$ and obtain a system of generators 
$$ G^{\gamma_\star}_a = gen^{\gamma_\star}_a\cup  
(\bigcup_{\gamma_0\leq\gamma <\gamma_\star} \tilde{B}^\gamma_a )     $$
for the ideal $\langle O^{\gamma_\star}_a \rangle$.

We now prove the existence of a normalized system of generators for ideal $\langle O^{\gamma_\star}_a \rangle$. Any germ 
$f\in  G^{\gamma_\star}_a \cap \Phi^{\gamma_\star}_a$ is of the form
$\tilde{f}=\tilde{u}^i_\alpha +\tilde{h}$ with $\tilde{h}\in\mathcal{F}^\gamma_a$ and $\gamma<\gamma^\star$.
According to Proposition 2.3 and the assumption step of induction, the germ $\tilde{h}$ is represented by
$$\tilde{h}=\tilde{q_1}\tilde{f}_{t_1} + \cdots + \tilde{q_p}\tilde{f}_{t_p}  ,   $$
where $\tilde{f}_{t_i}\in \tilde{B}^\gamma_a$, $\tilde{q_i}\in\mathcal{F}^\gamma_a$, and the germ $\tilde{r}\in\mathcal{F}^\gamma_a$ does not depend on principal variables of $\tilde{B}^\gamma_a$. 
 Then the germ $\tilde{f}^\star = \tilde{u}^i_\alpha + \tilde{r}$ is included in a set  $ben^{\gamma_\star}_a$. To do so with every germ in 
$G^{\gamma_\star}_a \cap \Phi^{\gamma_\star}_a$, we obtain an normalized system of generators
$$ \tilde{B}^{\gamma_\star}_a = ben^{\gamma_\star}_a\cup  
(\bigcup_{\gamma_0\leq\gamma <\gamma_\star} \tilde{B}^\gamma_a )     $$
for the ideal $\langle O^{\gamma_\star}_a \rangle$.

Let us take a point $c\sim\pi_{\gamma_\star}(a)$, then the  ideal $\langle O^{\gamma_\star}_a \rangle$  is isomorphic to the  ideal $\langle O^{\gamma_\star}_c \rangle$ of the algebra $\mathcal{F}^{\gamma_\star}_a$.
Indeed, if a function $f$ lies in $S$, then 
$\tilde{f}_a = \tilde{u}^i_\alpha|_a +\tilde{g}_a$ and $\tilde{f}_c = \tilde{u}^i_\alpha|_c +\tilde{g}_c$ because $S$ is a normalized set. 
Since $c\sim\pi_{\gamma_\star}(a)$, then $\tilde{g}_a=\tilde{g}_c$ and the ideal  $\langle O^{\gamma_\star}_a \rangle$  is isomorphic to the ideal
 $\langle O^{\gamma_\star}_c \rangle$. Therefore,  the ideal
  $\langle O^{\gamma_\star}_c \rangle$ has a normalized system of generators.

It is easy to show that a point $b$, such that $\pi_\gamma(b)=b_\gamma$ for all $\gamma\in\Gamma$, satisfies  (\ref{Db}) and the set 
$$ \tilde{B}_c = \bigcup_{\gamma_0\leq \gamma} \tilde{B}^\gamma_c $$
is a normalized system of generators for the differential ideal 
$\langle \langle \tilde{S} \rangle \rangle_c$ of $\mathcal{F}_c$. Therefore, the ideal  $\langle \langle \tilde{S} \rangle \rangle_c$ is soft.
By construction, we see that the set $\tilde{B}_c$ coincides with the orbit $O(lt \tilde{S}_c)$.  Thus $S$ is a passive system at $c\sim a$
and the theorem is proved.

\noindent
{\bf Theorem 5.2.} {\it Let $S=\{f_1,\dots,f_k\}\subset\mathcal{F}(V)$ be a differential system  with finite number of indeterminates. Suppose that $S$ is a normalized set and satisfies reducibility conditions (\ref{tau2}) on  $V$.
Then  
 the system $S$ is passive on $V$ and the set 
\begin{equation}\label{Man}          
\mathcal{M}=\{z\in V: f(z)=0, f\in O(S) \}   
\end{equation}
is a manifold in $\mathbb{J}$.}

Our proof is almost the same as the proof of Theorem 5.1. We employ the following denotation:
$$\gamma_0 = min\{\gamma\in\Gamma: O(S)\cap\Phi^\gamma(V)\neq 0\}, \quad O^\gamma = O(S)\cap \mathcal{F}^\gamma(V) , $$
$$ C^\gamma = O(S)\cap\Phi^\gamma(V) , \quad 
S^\gamma=O(S)\cap\mathcal{F}^\gamma(V) .$$
Let  $\langle O^\gamma \rangle$ be an ideal of the algebra  $\mathcal{F}^\gamma(V)$ generated by $O^\gamma$. 

 Using transfinite induction, we prove that for all $\gamma\geq\gamma_0$ there exists a normalized system of generators of the ideal $\langle O^\gamma \rangle$. Just as in the above theorem, we see that $O^{\gamma_0}$ is a normalized system of generators of the ideal $\langle O^{\gamma_0} \rangle$.
 
 Suppose that for each $\gamma_0\leq \gamma< \gamma_\star$ there exists a  normalized system  of generators $B^\gamma$ of the ideal $\langle O^{\gamma} \rangle$. We need to check the existences of such a system for $\gamma=\gamma_{\star}$.   At first one obtain a special system of generators $G^{\gamma_\star}$ of the ideal $\langle O^{\gamma_\star} \rangle$. 
  For this purpose, we consider the two cases again:
  
  (1) All leading terms of functions in $C^{\gamma_\star}$ are distinct.
   
   (2) There exist at least two function $f_i, f_j\in C^{\gamma_\star}$ such that $lt f_i = lt f_j$.
   
\noindent  
In the first case the set 
$$G^{\gamma_\star} = C^{\gamma_\star} \cup 
(\bigcup_{\gamma_0\leq\gamma <\gamma_\star} B^\gamma )     $$
is a system of generators of the ideal $\langle O^{\gamma_\star} \rangle$ 

In the second case there exist functions $f,g\in C^{\gamma_\star}$ such that 
$lt f =lt g$.  Thus there exist functions $f_p,f_q\in S$ and elements $\mu,\nu\in\mathbb{N}^n$ satisfying 
$$lt f = D^\mu (lt f_p) = D^\nu lt f_q = lt g .    $$  
It follows from condition of our theorem that the difference $f-g$ reduces  to the zero function modulo $S$.  One of the functions $f,g$ is included in a set
$gen^{\gamma_\star}$. In the same way we search for all pairs of functions in $C^{\gamma_\star}$ with equal leading terms, form the set $gen^{\gamma_\star}$ and obtain a system of generators 
$$ G^{\gamma_\star} = gen^{\gamma_\star}\cup  
(\bigcup_{\gamma_0\leq\gamma <\gamma_\star} B^\gamma )     $$
for the ideal $\langle O^{\gamma_\star} \rangle$. 

 We can then construct the set $B^{\gamma_\star}$ as follows. Any function
 $f\in G^{\gamma_\star}$ is the form $u^i_\alpha+h$, where $h\in\mathcal{F}^\gamma(V)$ with $\gamma<\gamma_\star$.
Furthermore the function $h$ is a polynomial in principal variables of $B^\gamma$ and coefficients of this  polynomial depend only on parametric variables.

Using Preposition 2.4, we write 
$$h = \sum q_i f_{t_i} + r ,   $$ 
where $f_{t_i}\in B^\gamma$, $q_i\in\mathcal{F}^\gamma(V)$, and the function $r\in \mathcal{F}^\gamma(V)$ depends only on parametric variables.
We include then the function $f^\star =u^i_\alpha + r$ in a set $ben^{\gamma_\star}$.
To do so with every function in 
$G^{\gamma_\star} \cap \Phi^{\gamma_\star}$, we obtain an normalized system of generators
$$ B^{\gamma_\star} = ben^{\gamma_\star}\cup  
(\bigcup_{\gamma_0\leq\gamma <\gamma_\star} B^\gamma )     $$
for the ideal $\langle O^{\gamma_\star} \rangle$.  
The set  
$$ B = \bigcup_{\gamma_0\leq \gamma} B^\gamma $$
is a normalized system of generators for the ideal 
$\langle \langle S \rangle \rangle$ of $\mathcal{F}(V)$
and this ideal is soft.
It is easy to see that the set $B$ coincides with the orbit $O(lt S)$.  Thus the system $S$ is a passive on $V$. From Proposition 2.7 it follows that the set (\ref{Man}) is a manifold in $\mathbb{J}$ 
and the theorem is proved.

\section{Examples}

We exhibit some examples assuming that $n=2$, $m=1$ and denoting by $u$ the variable $u_{00}$. The sets $U_n=\{u_{ij}:i+j=n\}$ form a partition of $U=\{u_{ij}\}_{i,j\in\mathbb{N}}$. We shall sometimes apply the usual terminology of differential equations. 

The smooth function 
\begin{equation}\label{sh}          
f = u_{11}-\sinh u   
\end{equation}
corresponds to the partial differential equation
\begin{equation}\label{eqsh}          
 u_{tx}-\sinh u = 0  .
\end{equation}
It is known (see \cite{Ibrag}) that vector fields 
$$X_1=  (u_{03}-\frac{1}{2}u_{01}^3)\frac{\partial}{\partial u}+\cdots  , \quad
X_2=(u_{05}-\frac{5}{2}u_{01}^2u_{03}-\frac{5}{2}u_{01}u_{02}^2+\frac{3}{8}u_{01}^5  )\frac{\partial}{\partial u}+\cdots$$
are higher symmetries of the equation (\ref{eqsh}).

Let $S_1$ be a differential system consisting of the functions $f$ and $h_1= u_{03}-\frac{1}{2}u_{01}^3$. 
We want to show that the ideal $I=\langle \langle S_1 \rangle \rangle$  is soft. For this purpose we shall construct a passive system generating the ideal $I$.
The functions $f$ and $h_1$ are orderly solvable with respect to $u_{11}$ and $u_{03}$ respectively. 
It is a straightforward calculation to check that the function $\tau(f,h_1)$ (given by (\ref{tau})) reduces to the function $f_1=u_{02}- \frac{1}{2}u_{01}^2\tanh(u)$ modulo $S_1$.  
Then an easy calculation shows that the function $\tau(f,f_1)$ reduces to the function $f_2=u_{10}- 2 \cosh(u)/u_{01}$ modulo $f$.  It is easy to see that the function $\tau({f_1,f_2})$ reduces to $0$ modulo the system $S=\{f_1,f_2\}$.
Furthermore the system $S$ generates the ideal $I$ and is passive.

We now find solutions of the system $S$. The function $f_1$ produces the ordinary differential equation 
$$ u_{xx} - \frac{1}{2}u_{x}^2\tanh(u) = 0    $$
having the first integral  $u_x/\cosh^2 u$. Using this integral and the equation
$$  u_t = 2\frac{\cosh u}{u_x} ,  $$
we obtain the implicit solution
$$\int \frac{du}{\sqrt{\cosh u}} = cx- 2t/c +c_1     $$
with $c,c_1\in\mathbb{R}$. 

Consider now the system $S_2$ consisting of the functions $f$ and $h_2= u_{05}-\frac{5}{2}u_{01}^2u_{03}-\frac{5}{2}u_{01}u_{02}^2+\frac{3}{8}u_{01}^5  $. A direct calculation shows  that the function $\tau({f,h_2})$ reduces to the function 
$$f_3= u_{04}-u_{01}u_{03}\tanh u+ \frac{1}{2}u_{02}^2\tanh u -\frac{3}{2}u_{01}^2u_{02} +\frac{3}{8}u_{01}^4\tanh u $$
  modulo $S_2$. Then the function $\tau(f,f_3)$ reduces to  
$$ f_4=u_{10} + \frac{4( u_{01}^3-2u_{03})\cosh u}
{8u_{01}u_{03}-4u_{02}^2-3u_{01}^4}  $$
  modulo $S_3=\{f_3,f_4\}$. It is possible to check that the function $\tau(f_3,f_4)$ reduces to $0$ modulo $S_4=\{f_3,f_4\} $, the system $S_4$ is passive and it generates the soft ideal $\langle \langle S_2 \rangle \rangle$.
  
The next example is closely connected with the equation (\ref{eqsh}) as well.
The set  
$$ v_{xxx} - \frac{2v_xv_{xx}}{v} +rv^4+s= 0, \qquad r,s\in\mathbb{R},  $$
is invariant manifold of the partial differential equation
$$ v_t = v_{xx}/v   $$
 as shown in \cite{KapBook}.  Using  the transformation   $v=\exp w$ we
 rewrite the last equations as
$$w_{xxx}+w_xw_{xx}-w^3_{x}+r\exp(3w)+s\exp(-w)=0, \qquad w_t +(\exp(-w))_{xx} = 0 .    $$
These equations correspond to two functions
$$f_5= u_{03}+u_{01}u_{02}-u_{01}^3+r\exp(3u)+s\exp(-u), 
\quad f_6=u_{10} + (u_{02}-u_{01}^2)\exp(-u) .     $$   
It is easy to check that the system $S_5=\{f_5,f_6\}$ generates a soft ideal  
$\langle \langle S_5 \rangle \rangle$ and is passive. The function $D_2(f_6)$ reduces to  the function $f_7=u_{11}+r\exp(2u)+s\exp(-2u)$ modulo $f_5$. This function lies in ideal $\langle \langle S_5 \rangle \rangle$ and corresponds to equation
  $$u_{tx}=\sinh(2u)   $$
with $r=-1/2$ and $s=1/2$.

\section{Conclusion} 

The above approach can be modified in different directions. As an example, we shall give some notations concerning the difference equations. 

We consider again the space $\mathbb{J}$ with the coordinate functions
$\{x_i\}_{i\in\mathbb{N}_n}$, $\{u^j_\alpha\}^{j\in\mathbb{N}_m}_{\alpha\in\mathbb{N}^n}$. 
Denote by  $\mathcal{C}^k(\mathbb{J})$ an algebra
of  $k$-times continuously differentiable functions on $\mathbb{J}$ which depend only on finitely many variables. 
  
Let $\sigma_1,\dots,\sigma_n$ be commuting endomorphisms of  $\mathcal{C}^k(\mathbb{J})$. Define an action of $\sigma_l$ on the coordinate functions by
$$\sigma_l(x_i) = \phi_{il}(x,u), \qquad  \sigma_l(u^j_\alpha)=u^j_{\alpha+e_l},     $$
where $l\in\mathbb{N}_n$, $x=(x_1,\dots,x_n)$, $u=(u^1,\dots,u^m)$, $e_1=(1,0,\dots,0),\dots,e_n=(0,\dots,1)$, $\phi_{il}\in \mathcal{C}^k(\mathbb{J})$. By definition, we have
$$ \sigma_l f(x_1,\dots,u^j_\alpha,\dots)=
f(\sigma_l(x_1),\dots,\sigma_l(u^j_\alpha),\dots)    $$
for any function $f\in\mathcal{C}^k(\mathbb{J})$.
For example, let $n=2, m=1$ and endomorphisms $\sigma_1$, $\sigma_2$ are given by the formulas:
$$\sigma_1(x_1)=x_1+h_1,\qquad  \sigma_1(x_2)=x_2,\qquad  
\sigma_2(x_1)=x_1,\qquad \sigma_2(x_2)=x_2+h_2 ,    $$
$$\sigma_1(u^1_\alpha)=u^1_{\alpha+e_1},\qquad\sigma_2(u^1_\alpha)=u^1_{\alpha+e_2} $$
with $\alpha=(\alpha_1,\alpha_2)\in\mathbb{N}^2$.


A subset $S$ of $\mathcal{C}^k(\mathbb{J})$ is called a system of difference  equations if every function $f\in S$ depends on at least one of variables $u^i_\alpha$. An ideal $I$ of the algebra $\mathcal{C}^k(\mathbb{J})$ is a $\sigma$-ideal if $\sigma_i(I)\subset I$ for all $i\in\mathbb{N}_n$. In particular, the set
$ \{\sigma^\alpha(f): f\in S, \alpha\in\mathbb{N}^n \}   $
generates a $\sigma$-ideal. 

It is possible to introduce soft $\sigma$-ideals and passive systems of difference equations. Dorodnistyn  recently found applications of Lie groups to finite-difference equations, meshes, and difference functionals \cite{Dorod}.

{\it Acknowledgements:} 
 This research was performed under the financial support of a grant from the Russian government for the conduct
 of research under the direction of leading scientists at the Siberian Federal University (Contract No. 14.U26.31.006)

\end{document}